\documentclass[12pt,a4paper]{article}

\usepackage{anysize}

\marginsize{3cm}{3cm}{3.5cm}{3.5cm}

\begin{document}

\title{Easy Proofs of Some Borwein Algorithms for $\pi$}

\author{Jes\'{u}s Guillera}

\date{}

\maketitle

\section{The Borweins' algorithms for $\pi$}
In 1987 Jonathan and Peter Borwein, inspired by the works of
Ramanujan, derived many efficient algorithms for computing $\pi$
(see \textbf{\cite{borweinagm}} and \textbf{\cite{borwein2}}).
However their proofs are difficult for the nonspecialist in the
theory of elliptic modular functions. We will see that by using only
a formula of Gauss's and elementary algebra we are able to prove the
correctness of two of them. The two algorithms we will consider are
the quadratic algorithm:
\[ d_0=1/\sqrt{2}, \quad r_0=1/2, \]
\[ d_{n+1}=\frac{1-\sqrt{1-d_n^2}}{1+\sqrt{1-d_n^2}}, \]
\[ r_{n+1}=(1+d_{n+1})^2 r_n -2^{n+1} d_{n+1}, \]
in which $r_n$ tends to $1/\pi$ and in each iteration the number of
exact digits is doubled, and the quartic algorithm:
\[ s_0=1/\sqrt[4]{2}, \quad t_0=1/2, \]
\[ s_{n+1}=\frac{1-\sqrt[4]{1-s_n^4}}{1+\sqrt[4]{1-s_n^4}}, \]
\[ t_{n+1}=(1+s_{n+1})^4 t_n -2^{2n+2} s_{n+1} (1+s_{n+1}+s_{n+1}^2), \]
in which $t_n$ tends to  $1/\pi$, and the number of correct digits
is quadrupled in each iteration.

\section{A Gauss formula}
At the age of 14, C. F. Gauss, starting with two numbers $a_0=a$ and
$b_0=b$, considered the iteration
\begin{equation}\label{agm}
a_{n+1}=\frac{a_n+b_n}{2}, \qquad b_{n+1}= \sqrt{a_n b_n},
\end{equation}
and defined the arithmetic-geometric mean of $a$ and $b$ as the
common limit
\begin{equation}\label{mab}
M(a,b)= \lim a_n = \lim b_n.
\end{equation}
In the year 1799, at the age of 22, by means of the elliptic
integrals
\[
\int_0^{\pi/2} \frac{d \theta}{ \sqrt{a_n^2 \cos^2 \theta+b_n^2
\sin^2 \theta} }, \qquad \int_0^{\pi/2} \frac{cos^2 \theta d
\theta}{ \sqrt{a_n^2 \cos^2 \theta+b_n^2 \sin^2 \theta} },
\]
he discovered the following wonderful formula \textbf{\cite[pp.
87-102]{arndt}}:
\begin{equation}\label{gauss}
\sum_{k=0}^{\infty} 2^{k} (a_k^2-b_k^2)=1-\frac{2M^2}{\pi},
\end{equation}
in which $a_0=1$, $b_0=1/\sqrt{2}$, $a_n$ and $b_n$ are given by
(\ref{agm}), and $M$ stands for the arithmetic-geometric mean of the
numbers $1$ and $1/\sqrt{2}$. Gauss's formula (\ref{gauss}) was
rediscovered in 1976 independently by E. Salamin and R. Brent (see
\textbf{\cite{brent}} and \textbf{\cite{salamin}}).

\section{Proof of the quadratic algorithm}
Our proof is based on the sequence
\begin{equation}\label{sucrsubn}
r_n=\frac{1}{\pi} \frac{M^2}{a_n^2}+\frac{1}{2 a_n^2}
\sum_{k=n}^{\infty} 2^k(a_k^2-b_k^2),
\end{equation}
where $a_0=1$, $b_0=1/\sqrt{2}$, $a_n$ and $b_n$ are defined as in
(\ref{agm}), and $M=M(1,1/\sqrt{2})$, that is, $M=\lim a_n=\lim
b_n$. It is easy to see that $\lim r_n=1/\pi$. Using relation
(\ref{gauss}) we also see that $r_0=1/2$. In addition we can easily
obtain
\[ r_n a_n^2-r_{n+1} a_{n+1}^2=2^{n-1} (a_n^2-b_n^2), \]
which leads to the recurrence
\begin{equation}\label{recrsubn}
r_{n+1}= \left( \frac{a_n}{a_{n+1}} \right)^2 r_n-2^{n+1} \left(
\frac{a_n}{a_{n+1}}-1 \right).
\end{equation}
To arrive at the Borweins' quadratic algorithm we write system
(\ref{agm}) in the form
\begin{equation}\label{agmsys}
\frac{b_n}{a_n}=2 \left(\frac{a_n}{a_{n+1}}\right)^{-1}-1, \qquad
\frac{b_{n+1}}{a_{n+1}}=\frac{a_n}{a_{n+1}} \sqrt{\frac{b_n}{a_n}}
\end{equation}
and deduce from it the relation
\begin{equation}\label{relab}
\frac{a_n}{a_{n+1}}-1=\sqrt{1-\frac{b_{n+1}^2}{a_{n+1}^2}}.
\end{equation}
Finally we can see that if we define
\begin{equation}\label{dsubn}
d_n=\sqrt{1-\frac{b_n^2}{a_n^2}},
\end{equation}
then
\begin{equation}\label{recud}
d_{n+1}=\frac{a_n}{a_{n+1}}-1=\frac{a_n-b_n}{a_n+b_n}=\frac{1-\sqrt{1-d_n^2}}{1+\sqrt{1-d_n^2}}
\end{equation}
and
\begin{equation}\label{recur}
r_{n+1}=(1+d_{n+1})^2 r_n -2^{n+1} d_{n+1}.
\end{equation}
Substituting the values of $a_0$ and $b_0$ in (\ref{dsubn}), we get
the initial value $d_0=1/\sqrt{2}$.

\section{Proof of the quartic algorithm}
From the relation (\ref{recur}) we have
\begin{equation}\label{r2nmas2}
r_{2n+2}=(1+d_{2n+2})^2 r_{2n+1} -2^{2n+2} d_{2n+2}
\end{equation}
and
\begin{equation}\label{r2nmas1}
r_{2n+1}=(1+d_{2n+1})^2 r_{2n} -2^{2n+1} d_{2n+1}.
\end{equation}
From (\ref{recud}) we get
\begin{equation}\label{invdsubn}
d_n=\sqrt{ 1-\left( \frac{1-d_{n+1}}{1+d_{n+1}} \right)^2 }=\frac{
2\sqrt{d_{n+1}} }{1+d_{n+1}}.
\end{equation}
If we define $s_n=\sqrt{d_{2n}}$ then $d_{2n+2}=s_{n+1}^2$ and
(\ref{invdsubn}) implies that
\begin{equation}\label{dsn}
d_{2n+1}=\frac{2 s_{n+1}}{1+s_{n+1}^2}.
\end{equation}
This allows us to rewrite (\ref{r2nmas2}) and (\ref{r2nmas1}) as
\begin{equation}\label{r2nmas2s}
r_{2n+2}=(1+s_{n+1}^2)^2 r_{2n+1} -2^{2n+2} s_{n+1}^2
\end{equation}
and
\begin{equation}\label{r2nmas1s}
r_{2n+1}=\left( 1+\frac{2s_{n+1}}{1+s_{n+1}^2} \right)^2 r_{2n}
-2^{2n+2} \frac{s_{n+1}}{1+s_{n+1}^2}.
\end{equation}
Combining (\ref{r2nmas2s}) and (\ref{r2nmas1s}) and defining
$t_n=r_{2n}$, we obtain
\begin{equation}\label{tsubn}
t_{n+1}=(1+s_{n+1})^4 t_n -2^{2n+2} s_{n+1} (1+s_{n+1}+s_{n+1}^2).
\end{equation}
On the other hand, from (\ref{recud}) we have
\begin{equation}\label{dsub2nplus1}
 d_{2n+1}=\frac{1- \sqrt {1-s_n^4}}{1+\sqrt{1-s_n^4}},
\end{equation}
which together with (\ref{dsn}) gives
\[ \sqrt[4]{1-s_n^4}=\frac {1-s_{n+1}}{1+s_{n+1}}. \]
Finally we obtain
\begin{equation}\label{recssubn}
s_{n+1}=\frac{1-\sqrt[4]{1-s_n^4}}{1+\sqrt[4]{1-s_n^4}}.
\end{equation}

\section{Other algorithms for $\pi$}
Jonathan and Peter Borwein have proved the following cubic analog
\cite{borwein3} of Gauss's formula (\ref{gauss}):
\begin{equation}\label{cubicgauss}
\sum_{k=0}^{\infty} 3^{k}(a_k^2-a_{k+1}^2)=\frac{1}{3}-\frac{M^2}{
\pi},
\end{equation}
where $a_0=1$, $b_0=1/2 \cdot \sqrt[3]{18-6\sqrt{3}}$,
\begin{equation}\label{itecub}
a_{n+1}=\frac{a_n+2b_n}{2}, \qquad b_{n+1}= \sqrt[3]{\frac{b_n(a_n^2
+a_n b_n+b_n^2)}{3}},
\end{equation}
and $M=\lim a_n=\lim b_n$. \\
\par Rewriting system (\ref{itecub}) in a way that is analogous to that used in (\ref{agmsys}) and
defining
\begin{equation}\label{defcub}
e_{n+1}=\frac{a_n}{a_{n+1}}, \qquad r_n=\frac{1}{\pi}
\frac{M^2}{a_n^2}+\frac{1}{a_n^2} \sum_{k=n}^{\infty}
3^{k}(a_k^2-a_{k+1}^2),
\end{equation}
we recover the following Borwein cubic algorithm for $\pi$:
\[ e_0=\sqrt{3}, \quad r_0=1/3, \]
\[ e_{n+1}=\frac{3}{1+\sqrt[3]{8-(e_n-1)^3}}, \]
\[ r_{n+1}=e_{n+1}^2 r_n -3^n (e_{n+1}^2-1), \]
in which $r_n$ tends cubically to  $1/\pi$. In \cite{vmoll} we also
find the following quartic analog of Gauss's formula (\ref{gauss}):
\begin{equation}\label{quargauss}
\sum_{k=0}^{\infty} 4^{k}(a_k^4-a_{k+1}^4)=\frac{1}{4}-\frac{3M^4}{4
\pi},
\end{equation}
where $ a_0=1$, $b_0=\sqrt[4]{12\sqrt{2}-16}$,
\begin{equation}\label{ite}
a_{n+1}=\frac{a_n+b_n}{2}, \qquad b_{n+1}=\sqrt[4]{\frac{a_n
b_n(a_n^2+b_n^2)}{2}},
\end{equation}
and $M=\lim a_n=\lim b_n$. \\
\par In the same way as before, using this quartic analog and defining
\begin{equation}\label{defcuar}
e_{n+1}=\frac{a_n}{a_{n+1}}, \qquad r_n=\frac{1}{\pi}
\frac{M^4}{a_n^4}+\frac{4}{3a_n^4} \sum_{k=n}^{\infty}
4^{k}(a_k^4-a_{k+1}^4),
\end{equation}
we can obtain the following algorithm for $\pi$:
\[ e_0=\sqrt{2}, \quad r_0=1/3, \]
\[ e_{n+1}=\frac{2}{1+\sqrt[4]{1-(e_n-1)^4}}, \]
\[ r_{n+1}=e_{n+1}^4 r_n -\frac{4^{n+1}}{3} (e_{n+1}^4-1), \]
in which $r_n$ tends quartically to  $1/\pi$. \\

{\it Av. Ces\'{a}reo Alierta, 31 esc. izda {\rm $4^o$}--A, Zaragoza
50008, Spain. \par jguillera@able.es}

\enddocument